\documentclass[]{article}
\usepackage{}
\usepackage{amsfonts}
\usepackage{latexsym,bm}
\usepackage{amssymb}
\usepackage{amsmath}
\usepackage{amsthm}
\usepackage{mathrsfs}
\usepackage[all]{xy}
\numberwithin{equation}{section}

\theoremstyle{plain}
\newtheorem{thm}{Theorem}[section]
\newtheorem{lem}[thm]{Lemma}
\newtheorem{prop}[thm]{Proposition}
\newtheorem{cor}[thm]{Corollary}
\theoremstyle{remark}
\newtheorem{rmk}[thm]{Remark}
\newtheorem{e.g.}[thm]{Example}
\theoremstyle{definition}
\newtheorem{defi}[thm]{Definition}
\newtheorem{conj}[thm]{Conjecture}

\begin{document}
\title{\Large Smooth projective horospherical varieties with nef tangent bundles}
\author{\normalsize Qifeng Li}
\maketitle

\begin{abstract}
We show that smooth projective horospherical varieties with nef tangent bundles are rational homogeneous spaces.
\end{abstract}



\section{\normalsize Introduction}

We work over the field of complex numbers. A famous conjecture of Hartshorne proved by Mori \cite{Mori79} in 1979 shows that projective spaces are the only smooth projective varieties with ample tangent bundles. A natural question is whether there is a similar characterization for smooth projective varieties with some semipositive conditions.  Demailly, Peternell and Schneider \cite{DPS94} proved that if $X$ is a compact K\"{a}hler manifold with nef tangent bundle, then there is a finite \'{e}tale cover $\widetilde{X}$ such that the Albanese map $\alpha: \widetilde{X}\rightarrow A(\widetilde{X})$ is a smooth fiberation with fibers being Fano manifolds with nef tangent bundles. The question above will be answered if the following conjecture of  Campana and Peternell \cite{CP91} holds.

\begin{conj} [Campana-Peternell Conjecture] \label{campana-peternell conjecture}
Smooth projective Fano varieties with nef tangent bundles are rational homogeneous spaces.
\end{conj}

From now on, we simply call a Fano manifold with nef tangent bundle a CP-manifold. Demailly, Peternell and Schneider \cite{DPS94} confirm Campana-Peternell Conjecture for Fano manifolds with dimension at most three. Through the works of Campana and Peternell \cite{CP93}, Mok \cite{Mok02}, Hwang \cite{Hw06}, Watanabe \cite{Wa14}\cite{Wa15}, Kanemitsu \cite{Ka15a}\cite{Ka15b}, Campana-Peternell Conjecture has been proved for Fano manifolds $X$ with $\dim(X)\leq \rho(X)+4$ , where $\rho(X)$ is the Picard number of $X$. The case when all the elementary Mori contractions of the manifold are smooth $\mathbb{P}^{1}$-fibrations have been solved by the series work of Mu\~{n}oz, Occhetta, Sol\'{a} Conde,  Watanabe and Wi\'{s}niewski in \cite{MOSW15} and \cite{OSWW14}. In particular, \cite{OSWW14} characterizes complete flag variety as the Fano manifolds all of whose elementary Mori contractions are smooth $\mathbb{P}^{1}$-fibrations.

In higher dimension very little is known except for the complete flag varieties cases. In fact, it is not known whether Campana-Peternell conjecture is true for quasi-homogeneous varieties. The aim of this paper is to verify the Campana-Peternell Conjecture for an important class of quasi-homogeneous varieties. We will show the following

\begin{thm} \label{thm. horospherical CP-mfds are homogeneous}
Smooth projective horospherical varieties with nef tangent bundles are rational homogeneous spaces.
\end{thm}

The sketch of the proof is as follows. In Section \ref{section Mori contractions on horospherical manifolds}, we study some special Mori contractions on horospherical varieties and reduce the proof of Theorem \ref{thm. horospherical CP-mfds are homogeneous} to the cases of Picard number one. In Section \ref{section Indices of smooth projective horospherical varieties of Picard number one}, we study the indices of smooth projective non-homogeneous horospherical varieties of Picard number one, and show that except for two cases they do not have nef tangent bundles. In Section \ref{section VMRT of a Mukai variety and odd symplectic Grassmannians}, we study the singularity of the VMRTs of the two exceptions, which imply that these varieties do not have nef tangent bundles.

\medskip

\textbf{\normalsize Notations and Conventions.} Denote by $G$ a connected reductive algebraic group. Let $B$ be a Borel subgroup of $G$, $T$ be a maximal torus in $B$, and $R_{u}(B)$ be the unipotent radical of $B$. Denote by $S$ the set of simple roots of $G$. For a subset $I$ of $S$, denote by $P_{I}$ the corresponding parabolic subgroup containing $B$ such that under that correspondence $P_{\emptyset}=B$ and $P_{S}=G$. When mentioning a rational homogeneous space $G/P_{I}$, we always mean its image under the minimal embedding unless otherwise stated.

When $G$ is a simple group, we use the standard notation for roots as in Bourbaki \cite{Bou75}. More precisely, $\omega_{i}$ is denoted for the $i$-th fundamental dominant weight, and $\alpha_{i}$ be the $i$-th simple root. We also write $P(\omega_{i})=P_{S\backslash\{\alpha_{i}\}}$.

For a linear algebraic group $H$, denote by $\mathscr{X}(H)$ the group of characters of $H$. For a subset $A$ of $\mathscr{X}(H)$, denote by $\text{Ker}_{H}(A)=\{h\in H \mid \chi(h)=0, \forall\chi\in A \}$. For an $H$-variety $X$ and a point $x\in X$, denote by $H_{x}$ the isotropy group of $x$ and $H\cdot x$ the $H$-orbit of $x$. Denote by $\mathbb{C}(X)^{(H)}=\{f\in \mathbb{C}(X)^{(H)}\backslash\{0\} \mid \exists \chi\in\mathscr{X} \text{\,s.t.\,} \forall h\in H, hv=\chi(h)f\}$, i.e. $\mathbb{C}(X)^{(H)}$ is the set of $H$-semiinvariant rational functions on $X$.

\section{\normalsize Mori contractions on horospherical manifolds} \label{section Mori contractions on horospherical manifolds}

\subsection{\normalsize Horospherical varieties}

In this subsection, we will review some properties on horospherical varieties. One can consult \cite{Pas08} for more details or  \cite{Per14} for related results on spherical varieties.

\begin{defi}
A $G$-variety $X$ is called $G$-horospherical if there exists a point $x_{0}\in X$ such that $G\cdot x_{0}$ is an open $G$-orbit on $X$ and the isotropy group $G_{x_{0}}\supseteq R_{u}(B)$. In this situation, we also say $X$ is a horospherical $G/H$-embedding, where $H=G_{x_{0}}$. Irreducible $B$-stable divisors having nonempty intersection with the open $G$-orbit are called colors of $X$. Denote by $\mathfrak{D}(G/H)$ the set of colors of $X$.
\end{defi}

Let $X$ be a $G$-horospherical variety. Let $x_{0}\in X$ and $H\subseteq G$ be as in the definition above. Then the normalizer $N_{G}(H)$ of $H$ is a parabolic subgroup of $G$ (see \cite[Propostion 2.2]{Pas08}). Denote by $I$ the subset of $S$ such that $N_{G}(H)=P_{I}$. From the Bruhat decomposition on $G/P_{I}$, we can get a natural bijective correspondence $S\backslash I\rightarrow\mathfrak{D}(G/H)$, $\alpha\mapsto D_{\alpha}$. When there is no confusions, we do not distinguish $S\backslash I$ and $\mathfrak{D}(G/H)$.

Denote by $M_{G/H}=\mathscr{X}(P_{I}/H)$. The character group $M_{G/H}$ is naturally a subgroup of $\mathscr{X}(B)$, and it is isomorphic to $\mathbb{C}(X)^{(B)}/\mathbb{C}^{*}$. Let $N_{G/H}=\text{Hom}(M_{G/H}, \mathbb{Z})$. For any $G$-orbit $Z$ on $X$, denote by $\mathfrak{C}_{Z}$ the associated cone in $(N_{G/H})_{\mathbb{Q}}$, and $\mathfrak{D}_{Z}=\{D\in\mathfrak{D}(G/H)\mid Z\subseteq D\}$. Denote by $\mathbb{F}_{X}$ the associated colored fan of $X$, i.e.
$\mathbb{F}_{X}=\{(\mathfrak{C}_{Z}, \mathfrak{D}_{Z})\mid Z$  is a $G$-orbit on $X\}$.
Denote by $\mathfrak{D}_{X}=\bigcup\mathfrak{D}_{Z}$, where $Z$ runs over the set of $G$-orbits on $X$.

The following lemma is collected from \cite[Proposition 2.4]{BM13} and \cite[Corollary 7.8]{Li15}.

\begin{lem} \label{lem. structure of orbits on horospherical varieties}
Let $X$ be a horospherical $G/H$-embedding. Assume $H\supseteq R_{u}(B)$. Denote by $I$ the subset of $S$ such that $P_{I}=N_{G}(H)$. Take a $G$-orbit $Z$ on  $X$. Then the following hold.

$(a)$ $M_{\mathfrak{C}_{Z}}\subseteq\mathscr{X}(P_{I\cup\mathfrak{D}_{Z}})$, where $M_{\mathfrak{C}_{Z}}=\{\chi\in M_{G/H} \mid \forall v\in\mathfrak{C}_{Z}, \langle \chi, v\rangle=0\}$.

$(b)$ $Z$ is $G$-equivariantly to $G/K$, where $K=\text{Ker}_{P_{I\cup\mathfrak{D}_{Z}}}M_{\mathfrak{C}_{Z}}\supseteq H\supseteq R_{u}(B)$.

$(c)$ $M_{\mathfrak{C}_{Z}}=M_{G/K}$ and $N_{G}(K)=P_{I\cup\mathfrak{D}_{Z}}$.

$(d)$ $\dim(Z)=\text{rank}(M_{\mathfrak{C}_{Z}})+\dim(G/P_{I\cup\mathfrak{D}_{Z}})$.
\end{lem}

Denote by
\begin{eqnarray*}
\mathfrak{D}_{0}(G/H)=\{D\in\mathfrak{D}(G/H)\mid \forall f\in\mathbb{C}(X)^{(B)}, \nu_{D}(f)=0\}.
\end{eqnarray*}

The following lemma is collected from Proposition 4.6, Theorem 4.7, Lemma 7.17, and Corollary 7.23 in \cite{Li15}. Note that for any parabolic subgroup $P$ of $G$ containing $B$, the intersection $P\cap P^{-}$ is a connected reductive algebraic group.

\begin{lem} \label{lem. structure of horospherical varieites}
Let $X$ be a smooth projective horospherical $G/H$-embedding, where $H\supseteq R_{u}(B)$ and $N_{G}(H)=P_{I}$.

$(i)$ Set $D_{0}=\sum\limits_{D\in\mathfrak{D}_{0}(G/H)}D$, and $J$ the subset of $S$ corresponding to $\mathfrak{D}_{0}(G/H)$. Then the linear system $|D_{0}|$ is base point free, and it induces a $G$-equivariant morphism $\pi: X\rightarrow G/P_{S\backslash J}$.

$(ii)$ The Mori cone $\overline{NE}(X)=F\times F'$, where $F$ is the extremal face corresponding to the Mori contraction $\pi$ and $F'$ is another extremal face.

$(iii)$ Set $L=P_{S\backslash J}\cap P_{S\backslash J}^{-}$. Then all fibers of $\pi$ are isomorphic to a same irreducible $L$-horospherical variety $Y$.

$(iv)$ If all effective divisors on $X$ are nef, then $S$ is the disjoint union of $I$, $J$ and $\mathfrak{D}_{X}$, and $Y$ is $L$-equivariantly isomorphic to $\prod\limits_{i=1}^{m}X_{i}$, where each $X_{i}$ is a smooth projective $L$-horospherical variety of Picard number one.
\end{lem}

\subsection{\normalsize Reduction to Picard number one cases}

The following results on manifolds with nef tangent bundles are basically due to \cite[Proposition 3.7, Theorem 5.2]{DPS94}. See also \cite[Theorem 4.4]{SW04} and \cite[Theorem 3.3]{MOSWW14} for modification of proofs of the conclusion $(ii)$.

\begin{prop} \label{prop. basic properties of CP-mfds}
Let $X$ be a CP-manifold. Then the following hold.

$(i)$ All effective divisors on $X$ are nef.

$(ii)$ Let $f: X\rightarrow Y$ be the contraction of an extremal ray. Then $Y$ and $f$ are smooth. Moreover, $Y$ and the fibers of $f$ are CP-manifolds.
\end{prop}

Now we can state and prove the main result in this section. This proposition reduces the proof of Theorem \ref{thm. horospherical CP-mfds are homogeneous} to the case when the horospherical CP-manifold is of Picard number one.

\begin{prop} \label{prop. horo. CP-mfd is product of homog. and those with rho=1}
Let $X$ be a smooth projective $G$-horospherical variety with nef tangent bundle. Then $X$ is isomorphic to $(\prod\limits_{i=1}^{m}X_{i})\times G/P$, where $P$ is a parabolic subgroup of $G$, and each $X_{i}$ is a smooth projective horospherical variety of Picard number one with nef tangent bundle.
\end{prop}

\begin{proof}
Step 1. In this step, we will construct two Mori contractions on $X$.

Horospherical varieties are rational. Hence, $X$ is simply connected and all finite \'{e}tale cover over $X$ is trivial (see \cite[Corollary 4.18]{De01}). By considering the Albanese map (see \cite[Theorem 3.14]{DPS94} or the statement of this result before Conjecture \ref{campana-peternell conjecture} in our introduction), $X$ is a CP-manifold. Hence, all effective divisors on $X$ are nef by Proposition \ref{prop. basic properties of CP-mfds}$(i)$. Thus, we can apply Lemma \ref{lem. structure of horospherical varieites}. There is a Mori contraction $\pi: X\rightarrow G/P_{S\backslash J}$ with fibers isomorphic to a same irreducible $L$-horospherical variety $\prod\limits_{i=1}^{m}X_{i}$, where $J=\mathfrak{D}(G/H)\subseteq S$, and each $X_{i}$ is a smooth projective $L$-horospherical variety of Picard number one.




The fact $X$ is a CP-manifold implies that $Y$ and all $X_{i}$ are CP-manifolds by Proposition \ref{prop. basic properties of CP-mfds}$(ii)$. So we only need to prove that $X\cong Y\times(G/P_{S\backslash J})$.


By Lemma \ref{lem. structure of horospherical varieites}$(ii)$, the Mori cone $\overline{NE}(X)=F\times F'$, where $F$ is the extremal face corresponding to the Mori contraction $\pi$ and $F'$ is another extremal face.

Denote by $\Psi: X\rightarrow W$ the Mori contraction associated with the extremal face $F'$. For the existence of such $\Psi$, see \cite[Theorem 3.1$(i)$]{Br93}. From the same theorem, we know that $W$ is a projective $G$-horospherical variety and $\Psi$ is a $G$-equivariant morphism. Moreover, $\Psi$ is a smooth fibration by Proposition \ref{prop. basic properties of CP-mfds}$(ii)$.  Then for any point $w\in\Psi(Z)$,
\begin{eqnarray} \label{equation dim of fiber=dX-dW}
\dim\Psi^{-1}(w)=\dim X-\dim W.
\end{eqnarray}

By the description of Mori contractions on horospherical varieties in terms of colored fans (see \cite[Proposition 3.4]{Br93}) and Lemma \ref{lem. structure of orbits on horospherical varieties}, we can get the conclusions $(a)(b)(c)$ as follows.

$(a)$ Take any $G$-orbit $Z$ on $X$, $\Psi^{-1}(\Psi(Z))=Z$, and $\mathfrak{D}_{\Psi(Z)}=\mathfrak{D}_{Z}$.

$(b)$ There exists horospherical subgroups $H'$ and $K'$ of $G$ such that $K'\supseteq H'\supseteq H$, and $Z$ and $\Psi(Z)$ are $G$-equiviriantly isomorphic to $G/H'$ and $G/K'$ respectively.

$(c)$ $N_{G}(H')=P_{I\cup\mathfrak{D}_{Z}}$, $N_{G}(K')=P_{I\cup J\cup\mathfrak{D}_{Z}}$, and $P_{I\cup\mathfrak{D}_{Z}}/H'\cong P_{I\cup J\cup\mathfrak{D}_{Z}}/K'$.

In summary, there is a commutative diagram as follows:
\begin{eqnarray*}
\xymatrix{
Z\ar@{=}[r]\ar[d]&G/H'\ar[r]\ar[d]&G/P_{I\cup\mathfrak{D}_{Z}}\ar[d] \\
\Psi(Z)\ar@{=}[r]&G/K'\ar[r]&G/P_{I\cup J\cup\mathfrak{D}_{Z}}.
}
\end{eqnarray*}


\medskip

Step 2. The aim of this step is to show a claim on the Dynkin diagram of $G$.

For any subset $A$ of $S$, denote by $\Gamma_{A}$ the full subdiagram of the Dynkin diagram $\Gamma$ of $G$ with vertices $A$.

\textit{Claim:} Each element in $S\backslash (I\cup J)$ and each element in $J$ lie in different connected components of $\Gamma$.

Take any $G$-orbit $Z$ on $X$. Then by the conclusions $(a)(b)(c)$ in Step 1, for any point $w\in\Psi(Z)$,
\begin{eqnarray} \label{equation dim of fiber=dP(IJZ)-dP(IZ)}
\dim\Psi^{-1}(w)=\dim P_{I\cup J\cup\mathfrak{D}_{Z}}-\dim P_{I\cup\mathfrak{D}_{Z}}.
\end{eqnarray}
By equations (\ref{equation dim of fiber=dX-dW}) and (\ref{equation dim of fiber=dP(IJZ)-dP(IZ)}),
\begin{eqnarray} \label{equation dim P(IJZ)-dim P(IZ)=dim X - dim W}
\dim P_{I\cup J\cup\mathfrak{D}_{Z}}-\dim P_{I\cup\mathfrak{D}_{Z}}=\dim X-\dim W.
\end{eqnarray}
Now apply formula (\ref{equation dim P(IJZ)-dim P(IZ)=dim X - dim W}) to the open $G$-orbit. Then for any $G$-orbit $Z$ on $X$,
\begin{eqnarray} \label{equation dim. of different fibers equal}
\dim(P_{I\cup J\cup\mathfrak{D}_{Z}})-\dim(P_{I\cup\mathfrak{D}_{Z}})=\dim(P_{I\cup J})-\dim(P_{I}).
\end{eqnarray}

Note that the union $I\cup J\cup\mathfrak{D}_{Z}$ is a disjoint union. By formula (\ref{equation dim. of different fibers equal}), for any $G$-orbit $Z$ on $X$, each element in $\mathfrak{D}_{Z}$ and each element in $J$ lie in different connected components of $\Gamma_{I\cup J\cup\mathfrak{D}_{Z}}$.

By Lemma \ref{lem. structure of horospherical varieites}$(iv)$, $S$ is the disjoint union of $I$, $J$ and $\mathfrak{D}_{X}$. Recall that $\mathfrak{D}_{X}=\bigcup\mathfrak{D}_{Z}$, where $Z$ runs over the set of $G$-orbits on $W$. Thus, the claim holds.

\medskip

Step 3. In this step, we will combine the two Mori contractions to get the isomorphism.

Recall that $\overline{NE}(X)=F\times F'$, and $\pi$ and $\Psi$ are the contractions corresponding to $F$ and $F'$ respectively. Hence, the natural morphism
\begin{eqnarray*}
& \Phi: X\rightarrow W\times (G/P_{S\backslash J}) \\
& x\mapsto (\Psi(x), \pi(x))
\end{eqnarray*}
is finite onto the image. On the other hand, by the Claim in Step 2,
\begin{eqnarray} \label{equation dim(G/PI)=dim(G/PS-J)=dim(PIJ)}
\dim(G/P_{I})=\dim(G/P_{S\backslash J})+\dim(G/P_{I\cup J}).
\end{eqnarray}
Denote by $G/K=\Psi(G/H)$ the open $G$-orbit on $W$ with $K\supseteq H$. Then by the conclusion $(c)$ in Step 1, $N_{G}(H)=P_{I}$, $N_{G}(K)=P_{I\cup J}$, $P_{I}/H\cong P_{I\cup J}/K$. Hence,
\begin{eqnarray*}
&\dim(X)&=\dim(G/P_{I})+\dim(P_{I}/H) \\
&&=\dim(G/P_{S\backslash J})+\dim(G/P_{I\cup J})+\dim(P_{I\cup J}/K) \\
&&=\dim(G/P_{S\backslash J})+\dim(W).
\end{eqnarray*}
Thus, the morphism $\Phi: X\rightarrow W\times (G/P_{S\backslash J})$ is surjective.

Take any point $x\in X$. Denote by $w=\Psi(x)\in W$. Then $\Psi^{-1}(w)$ is naturally a $G_{w}$-variety. By the conclusion $(a)$ in Step 1, $\Psi^{-1}(w)$ is $G_{w}$-homogeneous, and it is $G_{w}$-equivariantly isomorphic to $G_{w}/G_{x}$. By the conclusions $(a)(b)(c)$ in Step 1, $G_{w}/G_{x}\cong P_{I\cup J\cup\mathfrak{D}_{G\cdot w}}/P_{I\cup\mathfrak{D}_{G\cdot x}}$ and $\mathfrak{D}_{G\cdot w}=\mathfrak{D}_{G\cdot x}$. Thus, $\Psi^{-1}(w)$ is irreducible, it is isomorphic to $G_{w}^{o}/G_{x}\cap G_{w}^{o}$, and $G_{x}\cap G_{w}^{o}$ is a parabolic subgroup of the linear algebraic group $G_{w}^{o}$, where $G_{w}^{o}$ is the connected component of $G_{w}$ that contains the identity. Hence, the finite surjective $G_{w}^{o}$-equivariant morphism $\pi|_{\Psi^{-1}(w)}: \Psi^{-1}(w)\rightarrow G/P_{S\backslash J}$ is an isomorphism. So $\Phi: X\rightarrow W\times (G/P_{S\backslash J})$ is an isomorphism. By considering the Mori contraction $\pi$, we know $W$ is isomorphic to $Y$. Then the conclusion follows.
\end{proof}

\section{\normalsize Indices of smooth projective horospherical varieties of Picard number one} \label{section Indices of smooth projective horospherical varieties of Picard number one}

In this section, we turn to the cases of Picard number one. There is a classification of smooth projective horospherical varieties of Picard number one due to Pasquier.

\begin{prop} (\cite[Theorem 0.1, Lemma 1.19]{Pas09}) \label{prop. classification of horo. var. rho=1}
Let $X$ be a smooth projective $G$-horospherical variety of Picard number one. Assume that $X$ is not homogeneous. Then the following hold.

$(1)$ $G$ acts on $X$ with three orbits, one open orbit and two closed orbits. Identify them with $G/H$, $G/P_{1}$ and $G/P_{2}$ respectively, where $P_{1}$ and $P_{2}$ are parabolic subgroups of $G$ containing $B$. If necessary, reorder $P_{1}$ and $P_{2}$. Then the automorphism group $\text{Aut}(X)$ acts on $X$ with two orbits,  $G/H\cup G/P_{1}$ and $G/P_{2}$.

$(2)$ The variety $X$ is uniquely determined by the triple $(G, P_{1}, P_{2})$.  The triple $(G, P_{1}, P_{2})$ is one of the following list:

\quad\quad $(i)$ $(B_{m}, P(\omega_{m-1}), P(\omega_{m}))$ with $m\geq 3$;

\quad\quad $(ii)$ $(B_{3}, P(\omega_{1}), P(\omega_{3}))$;

\quad\quad $(iii)$ $(C_{m}, P(\omega_{k+1}), P(\omega_{k}))$ with $m\geq 2$ and $1\leq k\leq m-1$;

\quad\quad $(iv)$ $(F_{4}, P(\omega_{2}), P(\omega_{3}))$;

\quad\quad $(v)$ $(G_{2}, P(\omega_{2}), P(\omega_{1}))$,

\noindent where $G=B_{m}$ (resp. $C_{m}$, $F_{4}$, $G_{2}$) means $G$ is a simple group of that type.
\end{prop}

For smooth projective nonhomogeneous horospherical varieties of Picard number one, we have a formula as follows.

\begin{prop} \label{prop. r=dim(P1/P)+dim(P2/P)+2}
Keep notations $X, G, P_{1}, P_{2}$ as in Proposition \ref{prop. classification of horo. var. rho=1}. Denote by $r_{X}$ the Fano index of $X$, and set $P=P_{1}\cap P_{2}$. Then
\begin{eqnarray} \label{equation r=dim(P1/P)+dim(P2/P)+2}
r_{X}=\dim(P_{1}/P)+\dim(P_{2}/P)+2.
\end{eqnarray}
\end{prop}

\begin{proof}
Denote by $X_{i}$ the closed $G$-orbit on $X$ that is identified with $G/P_{i}$, and $U$ the open $G$-orbit on $X$. Consider the following diagram
\begin{eqnarray} \label{diagram blow-up X along X1 and X2}
\xymatrix{
\widetilde{X}\ar[d]_{\Phi}\ar[r]^{\pi}&G/P\\
X,
}
\end{eqnarray}
where $\Phi: \widetilde{X}\rightarrow X$ is the morphism of blowing-up $X$ along $X_{1}\cup X_{2}$, and the existence of the morphism $\pi$ follows from the fact that $X$ is a smooth horospherical variety. Moreover, both $\Phi$ and $\pi$ are $G$-equivariant morphisms, and $\pi$ is a $\mathbb{P}^{1}$-bundle. Denote by $E_{i}=\Phi^{-1}(X_{i})$. Then $\pi$ sends both $E_{1}$ and $E_{2}$ $G$-equivariantly isomorphically to $G/P$.

Set the colored fan of $X$ to be $\mathbb{F}_{X}=\{(U, \emptyset), (X_{1}, D_{2}), (X_{2}, D_{1})\}$, where $D_{1}$ and $D_{2}$ are the colors of $X$. Denote by $\widetilde{D}_{i}\subseteq\widetilde{X}$ the proper transform of $D_{i}$. Let $D_{ip}=\pi(\widetilde{D}_{i})$.

Now we will try to express the anticanonical divisor $-K_{G/P}$ using the information of $X$.

Since $X_{1}\nsubseteq D_{1}$ and $X_{2}\subseteq D_{1}$, we know that $\Phi^{*}(D_{1})=\widetilde{D}_{1}+a_{2}E_{2}$ in $\text{Pic}(\widetilde{X})_{\mathbb{Q}}$ for some $a_{2}\in\mathbb{Q}$. Now consider the following commutative diagram:
\begin{eqnarray} \label{diagram restriction on E_2}
\xymatrix{
E_{2}\ar[d]_{\Phi|_{E_{2}}}\ar[r]^{\pi|_{E_{2}}}&G/P\ar[d]^-{p_{2}}\\
X_{2}\ar[r]^-{\psi_{2}} &G/P_{2},
}
\end{eqnarray}
where $p_{2}$ is the natural morphism induced by the inclusion $P\subseteq P_{2}$, and $\psi_{2}$ is a $G$-equivariant isomorphism. Since  $\Phi: \widetilde{X}\rightarrow X$ is the blow-up morphism of $X$ along $X_{1}\cup X_{2}$, the restriction $\Phi|_{E_{2}}: E_{2}\rightarrow X_{2}$ is a $\mathbb{P}^{c_{2}-1}$-bundle, where $c_{2}$ is the codimension of $X_{2}$ in $X$. Thus,
\begin{eqnarray} \label{equation dim(P2/P)=c2-1}
\dim(P_{2}/P)=c_{2}-1,
\end{eqnarray}
and each fiber of $p_{2}$ is isomorphic to a projective space. Similarly, we know that
\begin{eqnarray} \label{equation dim(P1/P)=c1-1}
\dim(P_{1}/P)=c_{1}-1,
\end{eqnarray}
where $c_{1}$ is the codimension of $X_{1}$ in $X$.

By classical results on intersection theory on rational homogeneous spaces, we know that $\text{Pic}(G/P)=\mathbb{Z}(D_{1p})\oplus\mathbb{Z}(D_{2p})$, $NE(G/P)=\mathbb{Z}^{+}(l_{1})+\mathbb{Z}^{+}(l_{2})$, and $D_{ip}\cdot l_{j}=\delta_{ij}$. Moreover, $p_{2}$ is induced by the linear system $|D_{2p}|$. Hence, any line $l$ in an arbitrary fiber $F$ of $p_{2}$ is numerically equivalent to $l_{1}$. In particular, $D_{1p}\cdot l=1$ and $D_{1p}|_{F}$ is a generator of $\text{Pic}(F)$.

Denote by $\tilde{l}=(\pi|_{E_{2}})^{-1}(l)$. Since $\pi: \widetilde{X}\rightarrow G/P$ is a $\mathbb{P}^{1}$-bundle, $\widetilde{D}_{i}=\pi^{*}D_{ip}$ for $i=1, 2$. Thus, $\Phi^{*}(D_{1})\cdot\tilde{l}=(\widetilde{D}_{1}+a_{2}E_{2})\cdot\tilde{l}=D_{1p}\cdot l +a_{2}E_{2}\cdot\tilde{l}=1-a_{2}$.
Since $\Phi$ contracts $\tilde{l}$, $\Phi^{*}(D_{1})\cdot\tilde{l}=0$. Hence, $a_{2}=1$ and
\begin{eqnarray} \label{equation Phi^(-1)(D_1)=tilde(D)_1+E_2}
\Phi^{*}(D_{1})=\widetilde{D}_{1}+E_{2} \text{ in } \text{Pic}(\widetilde{X}).
\end{eqnarray}
Similarly, we have
\begin{eqnarray}
\Phi^{*}(D_{2})=\widetilde{D}_{2}+E_{1} \text{ in } \text{Pic}(\widetilde{X}).
\end{eqnarray}
Since $D_{1}=D_{2}$ in $\text{Pic}(X)$, we have
\begin{eqnarray}
\widetilde{D}_{1}-\widetilde{D}_{2}=E_{1}-E_{2} \text{ in } \text{Pic}(\widetilde{X}).
\end{eqnarray}
The fact $E_{1}\cap E_{2}=\emptyset$ implies that
\begin{eqnarray}
E_{1}\cdot E_{2} \text{ is a zero cycle on } \widetilde{X}.
\end{eqnarray}

\noindent The anticanonical divisor
\begin{eqnarray*}
-K_{X}=r_{X}D_{1}=r_{X}D_{2} \text{ in } \text{Pic}(X).
\end{eqnarray*}
Since $\Phi$ is the blow-up morphism, the anticanonical divisor
\begin{eqnarray}
-K_{\widetilde{X}}=\Phi^{*}(-K_{X})-(c_{1}-1)E_{1}-(c_{2}-1)E_{2} \text{ in } \text{Pic}(\widetilde{X}).
\end{eqnarray}
Moreover, the anticanonical divisor
\begin{eqnarray}
-K_{E_{1}}=(-K_{\widetilde{X}}-E_{1})|_{E_{1}} \text{ in } \text{Pic}(E_{1}).
\end{eqnarray}

Consider the following commutative diagram:
\begin{eqnarray*}
\xymatrix{
E_{1}\ar[rd]_{\pi|_{E_{1}}}\ar[r]^{j_{1}}&\widetilde{X}\ar[d]^{\pi}\\
&G/P,
}
\end{eqnarray*}
where $j_{1}$ is the natural inclusion morphism. Since the diagram is commutative and $\pi|_{E_{1}}$ is an isomorphism, we know that $(\pi|_{E_{1}})_{*}j_{1}^{*}\pi^{*}(D_{ip})=D_{ip}$ in $\text{Pic}(G/P)$. Recall that $\pi^{*}(D_{ip})=\widetilde{D}_{i}$. Hence,
\begin{eqnarray} \label{equation pi|_(E_1)(tilde(D)_i|_(E_1))=D_(ip)}
(\pi_{E_{1}})_{*}(\widetilde{D}_{i}|_{E_{1}})=D_{ip} \text{ in } \text{Pic}(G/P), \text{ for } i=1, 2.
\end{eqnarray}

Note that $\pi|_{E_{1}}: E_{1}\rightarrow G/P$ is an isomorphism. Combining with equations $(\ref{equation Phi^(-1)(D_1)=tilde(D)_1+E_2})-(\ref{equation pi|_(E_1)(tilde(D)_i|_(E_1))=D_(ip)})$, we can get that
\begin{eqnarray*}
-K_{G/P}=(r_{X}-c_{1})D_{1p}+c_{1}D_{2p} \text{ in } \text{Pic}(G/P).
\end{eqnarray*}
Similarly, by considering $E_{2}$ instead of $E_{1}$, we can get that
\begin{eqnarray*}
-K_{G/P}=c_{2}D_{1p}+(r_{X}-c_{2})D_{2p} \text{ in } \text{Pic}(G/P).
\end{eqnarray*}
Since $\text{Pic(G/P)}$ is freely generated by $D_{1p}$ and $D_{2p}$,  we know that
\begin{eqnarray*}
r_{X}=c_{1}+c_{2}=\dim(P_{1}/P)+\dim(P_{2}/P)+2,
\end{eqnarray*}
where the second equality follows from (\ref{equation dim(P2/P)=c2-1}) and (\ref{equation dim(P1/P)=c1-1}).
\end{proof}

\begin{prop} \label{prop. r>=r_i}
Keep notations as in Proposition \ref{prop. r=dim(P1/P)+dim(P2/P)+2}. Assume moreover that $X$ is a CP-manifold. Denote by $r_{i}$ the Fano index of $G/P_{i}$ for $i=1, 2$. Then
\begin{eqnarray} \label{equation r>=r_i}
r_{X}=\dim(P_{1}/P)+\dim(P_{2}/P)+2\geq\max\{r_{1}, r_{2}\}.
\end{eqnarray}
\end{prop}

\begin{proof}
Keep notation as in the proof of Proposition \ref{prop. r=dim(P1/P)+dim(P2/P)+2}. Let $C_{i}$ be a line on $X_{i}\cong G/P_{i}$, i.e. $C_{i}$ is an irreducible curve on $X_{i}$ with smallest anticanonical degree. Thus, $C_{i}$ is a smooth rational curve. Denote by $f_{i}: C_{i}\rightarrow X$ the natural inclusion. By assumption, the tangent bundle $\mathbb{T}_{X}$ is nef. Thus, $N_{X_{i}/X}|_{C_{i}}$, the restriction of the normal bundle, is nef. Consider the short exact sequence:
\begin{eqnarray*}
0\rightarrow \mathbb{T}_{X_{i}}|_{C_{i}}\rightarrow\mathbb{T}_{X}|_{C_{i}}\rightarrow N_{X_{i}/X}|_{C_{i}}\rightarrow 0.
\end{eqnarray*}
By the nefness of $N_{X_{i}/X}|_{C_{i}}$,
\begin{eqnarray*}
-K_{X}\cdot C_{i}=\deg(\mathbb{T}_{X}|_{C_{i}})\geq\deg(\mathbb{T}_{X_{i}}|_{C_{i}})=-K_{X_{i}}\cdot C_{i}=r_{i}.
\end{eqnarray*}

To complete the proof, we only need to show that $D_{i}\cdot C_{i}=1$ for $i=1,2$. This is equivalent to show that $D_{i}|_{X_{i}}$ is the ample generator of $\text{Pic}(X_{i})$, where the later is known to be true. For the convenience of the readers, we give a proof in detail as follows.

We consider the following commutative diagram:
\begin{eqnarray}
\xymatrix{
G/P\ar[d]_-{p_{2}}&E_{2}\ar[l]_-{\pi|_{E_{2}}}\ar[d]^-{\Phi|_{E_{2}}}\ar[r]^{j_{2}}&\widetilde{X}\ar[d]^-{\Phi}\\
G/P_{2}&X_{2}\ar[l]^-{\psi_{2}}\ar[r]_-{i_{2}} &X,
}
\end{eqnarray}
where $i_{2}$ and $j_{2}$ are natural inclusions.

Let $\widehat{D}_{2}=p_{2}(D_{2p})$. Then $\widehat{D}_{2}$ is the ample generator of $\text{Pic}(G/P_{2})$ and $p_{2}^{*}\widehat{D}_{2}=D_{2p}$. Thus,
\begin{eqnarray*}
&&(\Phi|_{E_{2}})^{*}i_{2}^{*}D_{2}=j_{2}^{*}\Phi^{*}D_{2}=(\widetilde{D}_{2}+E_{1})|_{E_{2}}=\widetilde{D}_{2}|_{E_{2}} \\ &&=(\pi|_{E_{2}})^{*}D_{2p}=(\pi|_{E_{2}})^{*}p_{2}^{*}\widehat{D}_{2}=(\Phi|_{E_{2}})^{*}\psi_{2}^{*}\widehat{D}_{2},
\end{eqnarray*}
where the fourth equality follows from an analogue of the formula (\ref{equation pi|_(E_1)(tilde(D)_i|_(E_1))=D_(ip)}).
Hence, $i_{2}^{*}D_{2}=\psi_{2}^{*}\widehat{D}_{2}$, which implies that $D_{2}|_{X_{2}}$ is the ample generator of $\text{Pic}(X_{2})$. Similarly, $D_{1}|_{X_{1}}$ is the ample generator of $\text{Pic}(X_{1})$. Then the conclusion follows.

\end{proof}

The inequality (\ref{equation r>=r_i}) can be checked for the five series of horospherical varieties $(i)-(v)$ in Proposition \ref{prop. classification of horo. var. rho=1}(2). Now keep notations as above. Denote by $d_{i}=\dim(P_{i}/P)$. Then we have the following table. It should be noticed that in Case $(i)$ of Proposition \ref{prop. classification of horo. var. rho=1}(2), we have $m\geq 3$.
\begin{center}
\begin{tabular}{|c|c|c|c|c|c|}
\hline
Case & $d_{1}$ & $d_{2}$ & $r_{1}$ & $r_{2}$ & Does (\ref{equation r>=r_i}) hold? \\ \hline
$(i)$ & $1$ & $m-1$ & $m+1$ & $2m$ & No \\ \hline
$(ii)$ & $3$ & $2$ & $5$ & $6$ & Yes \\ \hline
$(iii)$ & $2m-2k-1$ & $k$ & $2m-k+1$ & $2m-k$ & Yes \\ \hline
$(iv)$ & $2$ & $2$ & $5$ & $7$ & No \\ \hline
$(v)$ & $1$ & $1$ & $3$ & $5$ & No \\ \hline
\end{tabular}
\end{center}

As a direct consequence, we get the following

\begin{cor} \label{cor. horo. CP-mfd with rho=1 is either homog. or (ii) or (iii)}
Let $X$ be a CP-manifold as well as a horospherical variety of Picard one. Then $X$ is either a rational homogeneous space or a variety as in Proposition \ref{prop. classification of horo. var. rho=1}(2)$(ii)$ or $(iii)$.
\end{cor}


\begin{rmk} \label{rmk. geometric explanation of case (ii) (iii)}
varieties in Proposition \ref{prop. classification of horo. var. rho=1}$(2)(ii)$ and $(iii)$ have geometric explanations. Firstly, let $X$ be as in Case $(ii)$. By Proposition \ref{prop. r>=r_i} and the table above, $\dim(X)=9$ and $r_{X}=7$. Thus, $X$ is a Mukai variety. In fact, $X$ is isomorphic to a nonsingular section of the $10$-dimensional spinor variety $S^{10}\subseteq\mathbb{P}^{15}$ (see also \cite[Subsection 3.3]{FH12}).

Now assume $X$ to be as in Case $(iii)$. Then $X$ is isomorphic to the odd symplectic Grassmannian $G_{\omega}(k+1, 2m+1)$ with $m\geq 2$ and $1\leq k\leq m-1$ (see \cite[Proposition 1.12]{Pas09}). The odd symplectic Grassmannian is defined as follows. Let $V$ be a linear space of dimension $2m+1$. Equip a skew bilinear form $\omega$ of rank $2m$ on $V$. Then $G_{\omega}(i, 2m+1)$ is the variety of $i$-dimensional isotropic linear subspaces of $V$, where $1\leq i\leq m+1$.
\end{rmk}

\section{\normalsize Lines on horospherical varieties of Picard number one} \label{section VMRT of a Mukai variety and odd symplectic Grassmannians}

Let $X\subseteq\mathbb{P}^{N}$ be an $n$-dimensional (not necessarily irreducible or smooth) projective variety that are covered by lines. Denote by $F(X)$ the variety of lines in $X$. For any point $x\in X$, denote by $F(x, X)\subseteq\mathbb{P}(\mathbb{T}_{x}X)$ the variety of lines in $X$ passing through $x$, and we call it the variety of minimal rational tangents (VMRT for short) of $X$ at $x$. If moreover $F(x, X)\subseteq\mathbb{P}^{n-1}$ is covered by lines, then for any $[l]\in F(x; X)$, where $[l]$ stands for a line $l$ in $X$ passing through $x$, we denote $F(x, l, X)=F([l], F(x, X))$.

\begin{rmk}
$(i)$ The VMRT has the definition in a more general version. But we only need this restrictive version of definition here.

$(ii)$ Let $X$ be a smooth projective horospherical variety of Picard number one. Then the ample generator of $\text{Pic}(X)$, denoted by $\mathcal{O}_{X}(1)$, is very ample. Denote by $\Phi$ the closed embedding defined by the linear system $|\mathcal{O}_{X}(1)|$. Then the image $\Phi(X)$ is covered by lines. We always identify $X$ with the image $\Phi(X)$. In particular, the VMRTs are  naturally defined. Sometimes we also regard $X$ as an closed subvariety of a natural projective variety $Y$ (for example, in Proposition \ref{prop. VMRT of symplectic Grassmannian}). If we can show $F(x, X)\neq\emptyset$ for some $x\in X\subseteq Y$, then the inclusion $X\subseteq Y$ is compatible with the closed embedding induced by $|\mathcal{O}_{X}(1)|$.
\end{rmk}

The aim of this section is to show the following

\begin{prop} \label{prop. F(x, X) can be singular in both cases}
Let $X$ be a variety as in Case $(ii)$ or $(iii)$ in Proposition \ref{prop. classification of horo. var. rho=1}(2). Then there exists a point $x\in X$ such that $F(x, X)$ is singular.
\end{prop}

\begin{rmk} \label{rmk. all VMRTs of a CP-mfd are smooth}
If $X$ is a CP-manifold, then the nefness of the tangent bundle implies that any irreducible rational curve on $X$ is free. Hence, for any $x\in X$, $F(x, X)$ is smooth. Then by Proposition \ref{prop. F(x, X) can be singular in both cases}, varieties in Proposition \ref{prop. classification of horo. var. rho=1}(2)$(ii)$ and $(iii)$ do not have nef tangent bundles. This completes the proof of Theorem \ref{thm. horospherical CP-mfds are homogeneous}.
\end{rmk}

Let $V$ be an $n$-dimensional linear space. Denote by $\mathbb{F}(k_{1},\ldots,k_{m}; V)$ the flag variety parameterizing the sequences $(V_{k_{1}}, \ldots, V_{k_{m}})$ such that $V_{k_{1}}\subseteq V_{k_{2}}\subseteq\ldots\subseteq V_{k_{m}}$ and $V_{k_{i}}$ is a $k_{i}$-dimensional linear subspace of $V$. We also denote the Grassmannian $G(k, n)=G(k, V)=\mathbb{F}(k; V)$.

\begin{rmk} \label{rmk. VMRT of Grassmannian}
Keep notation as above. By \cite[Theorem 4.9]{LM03}, $F(G(k, V))=\mathbb{F}(k-1, k+1; V)$ and the VMRT of $G(k, V)\subseteq\mathbb{P}(\bigwedge^{k}V)$ is isomorphic to $\mathbb{P}^{k-1}\times\mathbb{P}^{n-k-1}\subseteq\mathbb{P}^{k(n-k)-1}$. Thus, the natural projection $\phi: \mathbb{F}(k-1, k, k+1; V)\rightarrow\mathbb{F}(k-1, k+1; V)$ is the universal family of lines on $G(k, V)$, and the natural projection $\psi: \mathbb{F}(k-1, k, k+1; V)\rightarrow G(k, V)$ is the evaluation morphism. Each line $l\subseteq G(k, n)$ corresponds to a point $x=(V_{k-1}, V_{k+1})\in\mathbb{F}(k-1, k+1; V)$, $l=\psi(\phi^{-1}(x))$ and $\phi^{-1}(x)=G(1, V_{k+1}/V_{k-1})$. So we can regard $l$ as a line in $G(1, V/V_{k-1})$. In particular, for any $k$-dimensional linear subspace $V_{k}$ of $V$, the VMRT $F([V_{k}], G(k, V))$ admits a morphism
\begin{eqnarray*}
&\Pi: F([V_{k}], G(k, V))\rightarrow G(k-1, V_{k})\\
 &[l] \mapsto [V_{k-1}],
\end{eqnarray*}
where $[l]=(V_{k-1}, V_{k+1})\in\mathbb{F}(k-1, k+1; V)$. The morphism $\Pi$ is a trivial $\mathbb{P}^{n-k-1}$- bundle and for each $(k-1)$-dimensional linear subspace $V_{k-1}$ of $V_{k}$, the fiber
\begin{eqnarray*}
\Pi^{-1}([V_{k-1}])=F([V_{k}/V_{k-1}], G(1, V/V_{k-1})).
\end{eqnarray*}
\end{rmk}

\medskip

Let $G=SO(7)$. Consider the following diagram:
\begin{eqnarray} \label{diagram G/P_(1 cap 3)}
\xymatrix{
G/(P(\omega_{1})\cap P(\omega_{3}))\ar[d]_{p_{3}}\ar[r]^-{p_{1}}&G/P(\omega_{1})\\
G/P(\omega_{3}).&
}
\end{eqnarray}
By the classical results on rational homogeneous spaces, $G/P(\omega_{1})$ and $G/P(\omega_{3})$ are smooth quadric hypersurfaces of dimension $5$ and $6$ respectively, $p_{1}$ is a $\mathbb{P}^{3}$-bundle, and $p_{3}$ is a $\mathbb{P}^{2}$-bundle. For any point $x\in G/P(\omega_{1})$, $p_{3}$ sends $p_{1}^{-1}(x)$ isomorphically to a $3$-dimensional linear subspace contained in the $6$-dimensional quadric hypersurface $G/P(\omega_{3})$.

\begin{lem} \label{lem. one dim. of V_3 have commen one dim. V_1}
Keep notations as above. Then there exists a line $l$ in $G/P(\omega_{3})$ and a line $C$ in $G/P(\omega_{1})$ such that for any $q\in C$, $p_{3}(p_{1}^{-1}(q))$ is a $3$-dimensional linear space containing $l$.
\end{lem}

\begin{proof}
Let $V$ be a $7$-dimensional linear space equipped with a nondegenerate symmetric bilinear form. Then for $1\leq i\leq 3$, the variety $G/P(\omega_{i})$ can be identified with the variety of $i$-dimensional isotropic linear subspaces of $V$. Take $V_{2}$ to be a $2$-dimensional isotropic linear subspace of $V$. Set
\begin{eqnarray*}
C:=G(1, V_{2}), \text{ and } l:= \{[V_{3}]\in G(3, V)\mid V_{2}\subseteq V_{3}\subseteq V_{3}^{\bot}\}.
\end{eqnarray*}
By regarding $G/P(\omega_{1})$ as a subvariety of the Grassmannian $G(1, V)$, we can see that $C$ is a line in $G/P(\omega_{1})$ (see Remark \ref{rmk. VMRT of Grassmannian}). Moreover, $l$ is a rational curve in $G/P(\omega_{3})$.

Now consider the following diagram:
\begin{eqnarray*}
\xymatrix{
G/(P(\omega_{2})\cap P(\omega_{3}))\ar[d]_{p_{3}}\ar[r]^-{p_{2}}&G/P(\omega_{2})\\
G/P(\omega_{3}).&
}
\end{eqnarray*}
Then $l=p_{3}(p_{2}^{-1}([V_{2}]))$. By the intersection theory on rational homogeneous spaces, $l$ is a line in $G/P(\omega_{3})$. Hence, for any $q\in C$, $p_{3}(p_{1}^{-1}(q))$ is a $3$-dimensional linear space containing $l$.
\end{proof}

\begin{prop}
Let $X\subseteq\mathbb{P}^{14}$ be a nonsingular hyperplane section of the $10$-dimensional spinor variety $S^{10}\subseteq\mathbb{P}^{15}$. Then there exists $x\in X$ such that $F(x, X)$ is not smooth.
\end{prop}

\begin{proof}
Firstly, $X$ is a horospherical variety as in Proposition \ref{prop. classification of horo. var. rho=1}(2)$(ii)$.
The automorphism group $\text{Aut}(X)$ acts on $X$ with two orbits. Denote by $Z$ the closed orbit. Then $Z\subseteq X$ is a $6$-dimensional smooth quadric hypersurface. Denote by $\Phi: \widetilde{X}\rightarrow X$ the blow-up of $X$ along $Z$. Let $E$ be the exceptional divisor of $\Phi$. Then there is a morphism $\pi: \widetilde{X}\rightarrow Q^{5}$ making $\widetilde{X}$ to be a $\mathbb{P}^{4}$-bundle over a $5$-dimensional smooth quadric hypersurface. In summary, we have the following commutative diagram:
\begin{eqnarray*}
\xymatrix{
E\ar@{^{(}->}[r]\ar[d]&\widetilde{X}\ar[r]^{\pi}\ar[d]_{\Phi}&Q^{5}\\
Z\ar@{^{(}->}[r]&X.&
}
\end{eqnarray*}
Moveover, the restriction $\pi|_{E}: E\rightarrow Q^{5}$ is a $\mathbb{P}^{3}$-bundle. Let $G=SO(7)$. Then $Z\cong G/P(\omega_{3})$, $G/P(\omega_{1})\cong Q^{5}$, $E\cong G/(P(\omega_{1})\cap P(\omega_{3}))$, and the morphisms $\pi|_{E}: E\rightarrow Q^{5}$ and $\Phi|_{E}: E\rightarrow Z$ coincide with $p_{1}$ and $p_{3}$ in the diagram (\ref{diagram G/P_(1 cap 3)}) respectively.

By Lemma \ref{lem. one dim. of V_3 have commen one dim. V_1}, there exists a line $C\subseteq Q^{5}$ and a line $l\subseteq Z$ such that for each $q\in C$, $\Phi((\pi|_{E})^{-1}(q))\subseteq Z$ is a $3$-dimensional linear subspace containing the line $l$.

Take a point $z\in l$. Since $X$ is a nonsingular hyperplane section of $S^{10}\subseteq\mathbb{P}^{15}$, $F(z, X)$ is a hyperplane section of $G(2, 5)\subseteq\mathbb{P}^{9}$. If $[l]$ is a singular point of $F(z, X)$, then $F(z, l, X)=\mathbb{P}^{1}\times\mathbb{P}^{2}\subseteq\mathbb{P}^{5}$. If $[l]$ is a smooth point of $F(z, X)$, then $F(z, l, X)$ is a hyperplane section of $\mathbb{P}^{1}\times\mathbb{P}^{2}\subseteq\mathbb{P}^{5}$.

Now take any point $q\in C$. Consider the following commutative diagram:
\begin{eqnarray*}
\xymatrix{
(\pi|_{E})^{-1}(q)\ar@{^{(}->}[r]\ar[d]&\pi^{-1}(q)\ar[d]\\
\Phi((\pi|_{E})^{-1}(q))\ar@{^{(}->}[r]&\Phi(\pi^{-1}(q)).
}
\end{eqnarray*}
By considering the diagram (\ref{diagram G/P_(1 cap 3)}), we get that $\Phi$ sends $(\pi|_{E})^{-1}(q)$ isomorphically to a 3-dimensional linear subspace containing the line $l$. Recall that $\pi^{-1}(q)$ is isomorphic to $\mathbb{P}^{4}$. Thus, $\Phi$ sends $\pi^{-1}(q)$ isomorphically to a 4-dimensional linear subspace containing the line $l$. Take any point $q'\in C\backslash\{q\}$, by the construction of the diagram (\ref{diagram G/P_(1 cap 3)}), $\Phi((\pi|_{E})^{-1}(q'))\neq\Phi((\pi|_{E})^{-1}(q))$. Thus, $\Phi(\pi^{-1}(q'))\neq\Phi(\pi^{-1}(q))$. Hence, there is a $1$-dimensional family of planes contained in $F(z, l, X)$. By the discussion in last paragraph, $F(z, l, X)=\mathbb{P}^{1}\times\mathbb{P}^{2}\subseteq\mathbb{P}^{5}$ and $F(z, X)$ is singular at the point $[l]$.
\end{proof}

\begin{prop} \label{prop. VMRT of symplectic Grassmannian}
Let $X$ be the odd symplectic Grassmannian $G_{\omega}(k, 2m+1)$ with $m\geq 2$ and $2\leq k\leq m$. Regard it as a subvariety of the Grassmannian $G(k, 2m+1)\subseteq\mathbb{P}(\bigwedge^{k}\mathbb{C}^{2m+1})$. Denote by $Z$ the closed orbit under the action of the automorphism group of $X$. Take a point $x\in X$. Then there is a surjective morphism $\pi: F(x, X)\rightarrow\mathbb{P}^{k-1}$.

$(i)$ If $x\notin Z$, then $\pi$ is a $\mathbb{P}^{2m-2k+1}$-bundle.

$(ii)$ If $x\in Z$, then there exists a hyperplane $H\subseteq\mathbb{P}^{k-1}$ such that
\begin{eqnarray*}
\pi^{-1}(v)\cong\left\{
\begin{array}{l@{\quad}l}
\mathbb{P}^{2m-2k+2}, & v\in H,\\
\mathbb{P}^{2m-2k+1}, & v\in\mathbb{P}^{k-1}\backslash H.
\end{array}\right.
\end{eqnarray*}
In particular, when $x\in Z$, $F(x, X)$ is not smooth.
\end{prop}

\begin{proof}
Let $V$ be a $(2m+1)$-dimensional linear space equipped with a skew bilinear form $\omega$ of rank $2m$. Denote by $M=G(k, 2m+1)\subseteq\mathbb{P}(\bigwedge^{k}V)$. Let $V_{k}$ be the $k$-dimensional linear subspace of $V$ represented by $x\in X\subseteq M$. By Remark \ref{rmk. VMRT of Grassmannian}, there exists a $\mathbb{P}^{2m-k}$-bundle $\Pi: F(x, M)\rightarrow G(k-1, V_{k})$ such that for each $(k-1)$-dimensional linear subspace $V_{k-1}$ of $V_{k}$,
\begin{eqnarray*}
\Pi^{-1}([V_{k-1}])=F([V_{k}/V_{k-1}], G(1, V/V_{k-1})).
\end{eqnarray*}
Restrict $\Pi$ on $F(x, X)$, we get the morphism $\pi: F(x, X)\rightarrow G(k-1, V_{k})=\mathbb{P}^{k-1}$.

For any linear subspace $W$ of $V$, define $W^{\bot}:=\{v\in V\mid \omega(v, W)=0\}$. Take an arbitrary $(k-1)$-dimensional linear subspace $V_{k-1}$ of $V_{k}$. The form $\omega$ induces a skew bilinear form $\tilde{\omega}$ on $V_{k-1}^{\bot}/V_{k-1}$. Moreover,
\begin{eqnarray*}
&&V_{k-1}^{\bot}\supseteq V_{k}^{\bot}\supseteq V_{k}\supseteq V_{k-1}, \\
&& X\cap G(1, V/V_{k-1})=G_{\tilde{w}}(1, V_{k-1}^{\bot}/V_{k-1}), \\
&&\pi^{-1}([V_{k-1}])=F([V_{k}/V_{k-1}], G_{\tilde{\omega}}(1, V_{k-1}^{\bot}/V_{k-1})).
\end{eqnarray*}

Note that $R:=V^{\bot}$ is a 1-dimensional linear subspace and the induced form on $V/R$ of the skew bilinear form $\omega$ is skew bilinear and nondegenerate. Hence,
\begin{eqnarray*}
\dim(V_{k-1}^{\bot})=\left\{
\begin{array}{l@{\quad}l}
2m-k+3, & R\subseteq V_{k-1},\\
2m-k+2, & R\nsubseteq V_{k-1},
\end{array}\right.
\end{eqnarray*}

and

\begin{eqnarray*}
G_{\tilde{\omega}}(1, V_{k-1}^{\bot}/V_{k-1}))=G(1, V_{k-1}^{\bot}/V_{k-1}))\cong\left\{
\begin{array}{l@{\quad}l}
\mathbb{P}^{2m-2k+3}, & R\subseteq V_{k-1},\\
\mathbb{P}^{2m-2k+2}, & R\nsubseteq V_{k-1}.
\end{array}\right.
\end{eqnarray*}

By \cite[Proposition 4.3]{Mih07}, $Z=\{[W]\in X\mid R\subseteq W\}$. Then the conclusion $(ii)$ holds, where $H$ is defined to be the hyperplane $G(k-2, V_{k}/R)$ in $G(k-1, V_{k})=\mathbb{P}^{k-1}$. Note that when $x\in Z$, $F(x, X)$ has two irreducible components with the same dimension and the intersection of them is not empty. So $F(x, X)$ is not smooth for $x\in Z$.

Now assume $x\notin Z$. Then $\pi$ is a $\mathbb{P}^{2m-2k+1}$-fibration over $\mathbb{P}^{k-1}$. It is known that such fiberation are projective bundles. So the conclusion $(i)$ holds.
\end{proof}

Finally, we summarize the proof of Theorem \ref{thm. horospherical CP-mfds are homogeneous} as follows.

\begin{proof}[Proof of Theorem \ref{thm. horospherical CP-mfds are homogeneous}]
Let $X$ be a smooth projective horospherical variety with nef tangent bundle. By Proposition \ref{prop. horo. CP-mfd is product of homog. and those with rho=1}, $X\cong(\prod\limits_{i=1}^{m}X_{i})\times G/P$ and each $X_{i}$ is a smooth projective horospherical variety of Picard number one as well as a CP-manifold.

It suffices to show that each $X_{i}$ is a rational homogeneous space. Now assume that there exists some $i_{0}$ such that $X_{i_{0}}$ is not homogeneous. By Corollary \ref{cor. horo. CP-mfd with rho=1 is either homog. or (ii) or (iii)}, $X_{i_{0}}$ is as in Proposition \ref{prop. classification of horo. var. rho=1}(2)$(ii)$ or $(iii)$. By Proposition \ref{prop. F(x, X) can be singular in both cases}, in both cases we can find $x\in X_{i_{0}}$ such that the variety $F(x, X_{i_{0}})$ is singular. This contradicts Remark \ref{rmk. all VMRTs of a CP-mfd are smooth}. Then the conclusion follows.
\end{proof}

\textbf{\normalsize Acknowledgements.} I want to thank Professor Baohua Fu and Professor Jun-Muk Hwang for communications related with this paper.

\small

Korea Institute for Advanced Study, 85 Hoegiro, Dongdaemun-gu, Seoul, 130-722, Republic of Korea

\smallskip

E-mail address: qifengli@kias.re.kr

\end{document}